\input amstex
\mag=\magstep1
\mag=\magstephalf % $B$3$l$O(J http://arxiv.org $B%5!<%P!<MQ(J
\documentstyle{amsppt}
\nologo
%\define\BaseLineSkip{\baselineskip=19.6pt}  
\define\BaseLineSkip{\baselineskip=13pt}  
\NoBlackBoxes
% \NoPageNumbers
% \NoRunningHeads
\rightheadtext\nofrills{DETERMINANT EXPRESSIONS FOR ABELIAN FUNCTIONS}
\leftheadtext\nofrills{YOSHIHIRO \^ONISHI}
\pagewidth{138mm}
%  \pageheight{205mm}
\hcorrection{2mm} % $B$3$l$O(J http://arxiv.org $BMQ(J
\vcorrection{-9mm} % $B$3$l$O(J http://arxiv.org $BMQ(J
% \leftskip=-7mm
\define\inbox#1{$\boxed{\text{#1}}$}
\def\fp{\flushpar}

\define\underbarl#1{\lower 1.4pt \hbox{\underbar{\raise 1.4pt \hbox{#1}}}}

\def\fp{\flushpar}
\define\tp#1{\negthinspace\ ^t#1}

\define\lr#1{^{\sssize\left(#1\right)}}

% ---- rightflush of \qed -------
\define\nullhbox{\hbox{\vrule height0pt depth 0pt width 1pt}}
\define\qedright{\null\nobreak\leaders\nullhbox\hskip10pt plus1filll\ \qed\par}
%\newbox\struthighbox
%\setbox\struthighbox=\hbox{\vrule height10pt depth 5pt width 0pt}
%\def\struthigh{\relax\ifmmode\copy\struthighbox\else\unhcopy\struthighbox\fi}

% ------ fonts setting -----------
%\jfont\tenptmbit=cmmib10 % Computer Modern Mathematical italic bold of 10pt
%\define\bk#1{\text{\tenptmbit{#1}}}
%\jfont\fourteengt=goth10 scaled \magstep5
%\jfont\twelveptgt=goth9 scaled \magstep2
%\jfont\tenptgt=goth8 scaled \magstep2
%\jfont\eightptgt=goth8
%\jfont\eightptmc=min8
%\jfont\sixptmc=min6 
\font\sc=cmcsc10
\font\twelveptrm=cmr12

\font\sc=cmcsc10
% ------ end fonts setting -------

\BaseLineSkip

\topmatter
\title  \nofrills\twelveptrm 
DETERMINANT EXPRESSIONS FOR ABELIAN FUNCTIONS IN GENUS TWO
\endtitle
\author {YOSHIHIRO \^ONISHI} \\
{\eightpoint \sl Faculty of Humanities and Social Sciences, Iwate University, 
                 Morioka, 020-8550, Japan} \\ 
{\eightpoint \sl e-mail:onishi\@iwate-u.ac.jp} 
\endauthor
% \address
% \endaddress
% \email
% \nofrills
% \centerline{\rm 
% \endemail
% \subjclass  
% \endsubjclass
% \abstract\nofrills{\bf

\endtopmatter
\document
\TagsOnRight
\BaseLineSkip

\par{\bf Abstract. }  
In this paper we generalize the formula of Frobenius-Stickelberger 
(see (0.1)  below) 
and the formula of Kiepert (see (0.2) below) to the genus-two case. 

$ $


$$ $$

\par{\bf Introduction. }  
There is a classical formula
  $$
%  \TagsOnRight
  \aligned
  &(-1)^{n(n-1)/2}1!2!\cdots n!
   \frac{\sigma(u\lr{0}+u\lr{1}+\cdots+u\lr{n})\prod_{i<j}\sigma(u_i-u_j)}
        {\sigma(u\lr{0})^{n+1}\sigma(u\lr{1})^{n+1}\cdots\sigma(u\lr{n})^{n+1}} \\
  &= \left|\matrix
      1 &  \wp(u\lr{0}) & \wp'(u\lr{0}) & \wp''(u\lr{0}) & \cdots & \wp\lr{n-1}(u\lr{0}) \\
      1 &  \wp(u\lr{1}) & \wp'(u\lr{1}) & \wp''(u\lr{1}) & \cdots & \wp\lr{n-1}(u\lr{1}) \\
 \vdots &  \vdots   & \vdots    & \vdots     & \ddots & \vdots           \\
      1 &  \wp(u\lr{n}) & \wp'(u\lr{n}) & \wp''(u\lr{n}) & \cdots & \wp\lr{n-1}(u\lr{n}) \\
   \endmatrix\right|
  \endaligned
  \tag 0.1
  $$
where  $\sigma(u)$  and  $\wp(u)$  are the usual elliptic functions.  
As far as the author knows the earliest work 
in which this formula appeared is a paper of 
Frobenius and Stickelberger \cite{\bf 10}.  
Before publication of their paper,  Kiepert \cite{\bf 14} 
gave the formula  
  $$
%  \aligned
  (-1)^{n-1}(1!2!\cdots (n-1)!)^2
   \frac{\sigma(nu)}{\sigma(u)^{n^2}}
  = \left|\matrix
    \wp'(u)        & \wp''(u)     &  \cdots  & \wp\lr{n-1}(u)  \\
    \wp''(u)       & \wp'''(u)    &  \cdots  & \wp\lr{n}(u)    \\
    \vdots         & \vdots       &  \ddots  & \vdots          \\
    \wp\lr{n-1}(u) & \wp\lr{n}(u) &  \cdots  & \wp\lr{2n-3}(u) \\
   \endmatrix\right|,   
%  \endaligned
  \tag 0.2
  $$
which also can be obtained by a limiting process from (0.1). 
These formulae are included in the book \cite{{\bf 21}, p.458 and p.460} 
as exercises and are also treated in \cite{{\bf 9}, I, p.183 and p.186}.  
The function  $\sigma(nu)/\sigma(u)^{n^2}$  in  (0.2) has been denoted 
by  $\psi_n(u)$  ever since the book of Weber \cite{\bf 20}. 

There is also a detailed analytic theory of Abelian functions 
for the Jacobian variety of any hyperelliptic curve.  
In the theory there is a nice generalization by H.F. Baker and others 
of the elliptic sigma function.  
If the genus of the curve, say  $C$, is  $g$, this function is also denoted 
by  $\sigma(u)=\sigma(u_1,\cdots, u_g)$ (see (1.1) below).  
Although the function
  $$
  \frac{\sigma(nu)}{\sigma(u)^{n^2}} 
  $$
is a function on the Jacobian variety and 
is a natural generalization of  $\psi_n(u)$  for the given genus  $g$  curve, 
it has poles along a theta divisor  $\Theta$  
which is the  $(g-1)$-fold sum of the curve  $C$  
embedded with a Weierstrass point as base point,  
so does not naturally restrict to a function on the curve itself.  

Now we restrict ourselves to consider only functions 
associated to a hyperelliptic curve of genus  $g=2$.  
Then we consider the modified function 
  $$
  \frac{\sigma(nu)}{\sigma_2(u)^{n^2}}, 
  $$
where  $\sigma_2(u)=(\partial\sigma/\partial u_2)(u)$.  
Although this function is not a function on the Jacobian variety 
of the curve, via restriction to  $C$  
it can be regarded as a function on the curve itself,  
and is a good candidate to be called the {\it hyperelliptic psi function}.  
So we denote it by  $\psi_n(u)$.  
The zeroes of  $\psi_n(u)$  on the curve are just the points 
whose multiplication by  $n$  are on  $\Theta$. 
If  $g=2$  the theta divisor  $\Theta$  coincides exactly with  $C$.   

The aim of this paper is to give a natural generalization of 
the expressions (0.1) and (0.2) for any hyperelliptic curve of genus two.  
See Theorem 2.3 and Theorem 3.3 below.  

The function  $\psi_n(u)$  has been investigated in several contexts 
as follows.  
First, it had an important role in the new theory of 
complex multiplication due to D. Grant (\cite{\bf 12} and \cite{\bf 18}). 
Secondly, a recursion relation for any hyperelliptic 
psi function was given by D.G. Cantor in \cite{{\bf 7}}, 
and such relation is used to compute the torsion points on 
the corresponding curve.  
He also gave a determinant expression  for  $\psi_n(u)$, 
which is different from our expression and 
should be regarded as a generalization of 
the formula of Brioschi \cite{{\bf 6}, p.770, $\ell$.3} 
to all hyperelliptic curves.

The formula (0.1) also has a classical generalization to 
all non-singular algebraic curves in terms of the Klein prime form, 
which is the formula (44) in \cite{{\bf 8}, p.33}, 
and a generalization to all hyperelliptic curves of genus two in terms of 
Gunning's prime form, which is the formula (2) of \cite{{\bf 13}}.  
The referee of this paper kindly informed the author 
about these two kinds of generalizations. 

The author would like to acknowledge that 
the proof of the main theorem (Theorem 2.3) was quite simplified 
by the referee.  

The first thing which inspired the author in this work are 
the very important and strange formulae 
which appeared in the famous papers \cite{\bf 15} and \cite{\bf 16}.  
These papers solved the problem of the determination of 
the argument of cubic and quartic Gauss sums, respectively,  
by using the formulae given by Propositions 7.3 and 7.4 
in \cite{{\bf 15}, p.181} and by (6.37) and (6.38) in \cite{\bf 16}, 
which resemble quite closely (0.2) above.  
Matthews mentions in \cite{{\bf 15}, p.179} that 
such a line of attack was suggested by observations of S.J. Patterson.  
The author hopes that our new formulae will be used to establish 
a new theory of complex multiplication.  

The author was led to believe in the existence of our determinant 
expressions through continual communication with Matsutani.  
The author found another generalization 
of the formula of Frobenius and Stickelberger 
to the three variables case (i.e. $n=2$) for genus-two curves 
in \cite{{\bf 5}, pp.96-97} after having already obtained 
these results.  

{\eightpoint {\it Added for the finally revised version}. 
Quite recently all our formulae were generalized to any hyperelliptic curves 
by the author himself in \cite{{\bf 19}}. 
Moreover the connection of Theorem 3.3 below and 
the determinant expression of Cantor mentioned above 
was completely clarified by S. Matsutani 
(see Appendix of \cite{{\bf 19}}).  
 }

$ $

\par{\bf 0. Convention. }
We denote, as usual, by  $\bold Z$ 
% , $\bold Q$  
and  $\bold C$ 
the ring of rational integers
%, the field of rational numbers 
and the field of complex numbers, respectively.  
In an expression of the Laurent expansion of a function, 
the symbol  $(d^{\circ}(z_1, z_2, \cdots, z_m)\geq n)$  
stands for the terms of total degree 
at least  $n$  with respect to 
the given variables  $z_1$, $z_2$, $\cdots$, $z_m$. 
% When the variable or the least total degree are clear from the context, 
% we simply denote them by  $(d^{\circ}\geq n)$  or the dots \lq\lq $\cdots$". 

In cross references, we indicate a formula as (1.2),  
and each of Lemmas, Propositions, Theorems and Remarks also as 1.2.  

$ $
    
\par{\bf 1. The Sigma Function in Genus Two. }  
In this Section we summarize the fundamental facts used in Sections 2 and 3.  
Detailed treatment of these facts are given 
in \cite{\bf 1}, \cite{\bf 2} and \cite{\bf 3} 
(see also Section 1 of \cite{\bf 18}).  
   
Let
  $$
  f(x)=\lambda_0     + \lambda_1 x   + \lambda_2 x^2 
     + \lambda_3 x^3 + \lambda_4 x^4 + \lambda_5 x^5
  $$
be an polynomial of  $x$  over  $\bold C$  such that its all zeroes are 
different each other.  
Let  $C$  be a smooth projective model of the curve of genus  $2$ 
defined by  $y^2=f(x)$.  
We denote by  $\infty$  the unique point at infinity. 
In this paper we suppose  $\lambda_5=1$.  The set of forms
  $$
  \omega_1=\frac{dx}{2y}, \ \omega_2=\frac{xdx}{2y}
  $$
is a basis of the space of differential forms of first kind.  
Let  
  $$
  \eta_1=\frac{(\lambda_3 x +2\lambda_4 x^2 + 3\lambda_5 x^3)dx}{2y}, \ \ 
  \eta_2=\frac{           x^2 dx}{2y}.
  $$
Then  $\eta\lr{1}$  and  $\eta\lr{2}$  are differential forms 
of the second kind without poles except at  $\infty$.  
% (\cite{B-abel, p.195, Ex.i} or \cite{B-sigma, p.314}).  
We fix generators  $\alpha_1$,  $\alpha_2$,  $\beta_1$, $\beta_2$  
of the fundamental group of  $C$  such that their intersections 
are  $\alpha_i\cdot\alpha_j=\beta_i\cdot\beta_j=0$,  
$\alpha_i\cdot\beta_j=\delta_{ij}$  for  $i$, $j=1$, $2$.  
If we set
  $$
  \omega' =\left[\matrix\ \int_{\alpha_1}\omega_1 & 
                          \int_{\alpha_2}\omega_1 \\
                          \int_{\alpha_1}\omega_2 &
                          \int_{\alpha_2}\omega_2
            \endmatrix\right], \ 
  \omega''=\left[\matrix\ \int_{\beta_1}\omega_1  &
                          \int_{\beta_2}\omega_1  \\
                          \int_{\beta_1}\omega_2  & 
                          \int_{\beta_2}\omega_2
            \endmatrix\right] 
  $$
the lattice of periods of our Abelian functions appearing below is given by 
  $$
  \Lambda=\omega' \left[\matrix \bold Z \\ \bold Z \endmatrix\right]
          +\omega''\left[\matrix \bold Z \\ \bold Z \endmatrix\right]
  (\subset \bold C^2). 
  $$
The modulus of  $C$  is  $Z:={\omega'}^{-1}\omega''$. 
We also introduce matrices 
  $$
  \eta' =\left[\matrix\ \int_{\alpha_1}\eta_1  &
                        \int_{\alpha_2}\eta_1 \\
                        \int_{\alpha_1}\eta_2  &
                        \int_{\alpha_2}\eta_2
           \endmatrix\right], \ 
  \eta''=\left[\matrix\ \int_{\beta_1} \eta_1  &
                        \int_{\beta_2} \eta_1 \\
                        \int_{\beta_1} \eta_2  &
                        \int_{\beta_2} \eta_2
           \endmatrix\right].
  $$

   Let  $J$  be the Jacobian variety of the curve  $C$.  
We identify  $J$  with the Picard group  $\text{Pic}^{\circ}(C)$  of 
linear equivalence classes of divisors of degree 0 of  $C$.  
Let  $\text{Sym}^2(C)$  be the symmetric product of two copies of  $C$.  
Then we have a birational map
  $$
  \align
  \text{Sym}^2(C)&\rightarrow \text{Pic}^{\circ}(C)=J \\
  (P_1, P_2) &\mapsto \text{the class of}\  P_1+P_2-2\cdot\infty.
  \endalign
  $$
We may also identify (the  $\bold C$-rational points of)  $J$  
with  $\bold C^2/\Lambda$.  
We denote by  $\kappa$  the canonical 
map  $\bold C^2\rightarrow \bold C^2/\Lambda$  
and by  $\iota$  the embedding of  $C$  into  $J$  given by mapping  $P$  
to  $\text{the class of}\ P-\infty$.  
Then the image  $\iota(C)$  is a theta divisor 
which is usually denoted by  $\Theta$. 
Although we treat only the case of genus 2, we use 
both of these two symbols to make as clear as 
possible the shape of our expected generalization to 
the case of higher genus.  
We denote by  $O$  the origin of  $J$.  
Obviously  $\Lambda=\kappa^{-1}(O)=\kappa^{-1}\iota(\infty)$.  
If  $u\in \bold C^2$, we denote by  $u_1$  and  $u_2$ 
the first and second entries of  $u$.

\proclaim{\indent\sc Lemma 1.1} 
As a subvariety of  $J$,  the divisor  $\Theta$  is non-singular. 
\endproclaim

\fp A proof of this fact is seen, for instance, 
in Lemma 1.7.2 of \cite{{\bf 18}, p.390}.  
If we set 
  $$
  \delta''=\left[\frac 12 \ \ \frac 12 \right], \ \ 
  \delta' =\left[    0    \ \ \frac 12 \right], \
%  \text{and}\ 
%  \delta:=\left[\matrix \delta'' \\ \delta' \endmatrix\right]
  $$
then the sigma function attached to  $C$  is defined, as in \cite{\bf 2}, by 
  $$
  \aligned
  \sigma&(u)=c\exp(-\frac 12u\eta'{\omega'}^{-1}\tp{u}) \\
  &\cdot\sum_{n\in \bold Z^2}
  \exp[2\pi \sqrt{-1}\{\frac 12\tp{(n+\delta'')}Z(n+\delta'')
                      +\tp{(n+\delta'')}({\omega'}^{-1}\tp{u}+\delta')\}].
  \endaligned
  \tag 1.1
  $$
To fix the constant  $c$  above we need the following lemma.

\proclaim{\indent\sc Lemma 1.2}
The Taylor expansion of  $\sigma(u)$  at  $u=(0,0)$  is of the form 
  $$
  u_1+\frac 16 \lambda_2 {u_1}^3 
  - \frac 13 \lambda_5 {u_2}^3 + (d^{\circ}(u_1, u_2)\geq 5)
  $$
up to a multiplicative constant.  
\endproclaim

\fp
Lemma 1.2 is proved in \cite{{\bf 4}, p.96} 
(see also \cite{{\bf 11}, pp.129-130} 
or Proposition 2.1.1(2) in \cite{{\bf 18}}).  
We fix the constant  $c$  in (1.1) such that the expansion is 
exactly of the form in 1.2.

\proclaim{\indent\sc Lemma 1.3} 
Let  $\ell$  be an element in  $\Lambda$.  
The function  $u \mapsto \sigma(u)$  on  $\bold C^2$  satisfies 
the translational formula
  $$
  \sigma(u+\ell)=\chi(\ell)\sigma(u)\exp L(u+\ell, \ell), 
  $$
where  $\chi(\ell)=\pm 1$  is independent of  $u$,   
$L(u,\  v)$  is a form which is bilinear over the field of real numbers 
and  $\bold C$-linear with respect to the first variable  $u$  
and  $L(\ell\lr{1},\ \ell\lr{2})$  is  $2\pi\sqrt{-1}$  times an integer 
if  $\ell\lr{1}$  and  $\ell\lr{2}$  are in  $\Lambda$.  \fp
\endproclaim

\fp The details of 1.3 are given in  \cite{{\bf 1}, p.286} 
(see also \cite{{\bf 18}, pp.395-396}).  

\proclaim{\indent\sc Lemma 1.4} 
{\rm (1)}  The function  $\sigma(u)$  on  $\bold C^2$  
vanishes if and only if  $u\in \kappa^{-1}(\Theta)$.  \fp 
{\rm (2)}  Suppose  $v\lr{1}$, $v\lr{2}$  are two points of  $\kappa^{-1}\iota(C)$.  
The function  $u\mapsto \sigma(u-v\lr{1}-v\lr{2})$  for  $u\in\kappa^{-1}\iota(C)$  
is identically zero if and only if  $v\lr{1}+v\lr{2}$  
is contained in  $\kappa^{-1}(O)(=\Lambda)$.  
If the function is not identically zero, 
it vanishes only at  $u=v\lr{1}$  and  $v\lr{2}$  modulo  $\Lambda$  to order  $1$  
or to order $2$ if these two points coincide.  \fp
{\rm (3)}  Let  $v$  be a fixed point of  $\kappa^{-1}\iota(C)$.  
There exist a point  $v\lr{1}$  of  $\kappa^{-1}\iota(C)$  
such that the function   $u\mapsto \sigma(u-v-v\lr{1})$  
on  $\kappa^{-1}\iota(C)$  is 
not identically zero and vanishes at  $u=v$  modulo  $\Lambda$  of order  $1$. 
\endproclaim

\fp 
The assertions 1.4(1) and (2) are proved in \cite{{\bf 1}}, pp.252-258, 
for instance. 
The assertion 1.4(3) obviously follows from (2).  

We introduce the functions 
  $$
  \wp_{jk}(u)=-\frac{\partial^2}{\partial u_j\partial u_k}
               \log\sigma(u), \ 
  \wp_{jk\cdots r}(u)=\frac{\partial}{\partial u_j}\wp_{k\cdots r}(u). 
  $$
These functions
are periodic with respect to the lattice  $\Lambda$ by 1.3, 
and has poles along  $\kappa^{-1}(\Theta)$  by 1.4(1).  
We also use the notation
  $$
  \sigma_j(u)=\frac{\partial}{\partial u_j}\sigma(u), \ 
  \sigma_{jk}(u)
         =\frac{\partial^2}{\partial u_j\partial u_k}\sigma(u).  
  $$
The following formula was also obtained by Baker in \cite{\bf 2}, p.381. 

\proclaim{\indent\sc Lemma 1.5} 
We have that
  $$
  -\frac{\sigma(u+v)\sigma(u-v)}{\sigma(u)^2\sigma(v)^2}
  =\wp_{11}(u)-\wp_{11}(v)+\wp_{12}(u)\wp_{22}(v)-\wp_{12}(v)\wp_{22}(u).
  $$
\endproclaim

Let  $(u_1, u_2)$  be an arbitrary point 
not in  $\kappa^{-1}(\Theta)$.      
Then we can find a unique pair of points  $(x_1, y_1)$  and  $(x_2, y_2)$  
on  $C$  such that 
  $$
  u_1=\int_{\infty}^{(x_1, y_1)}\omega_1 
        + \int_{\infty}^{(x_2, y_2)}\omega_1, \ \ 
  u_2=\int_{\infty}^{(x_1, y_1)}\omega_2 
        + \int_{\infty}^{(x_2, y_2)}\omega_2
  \tag 1.2
  $$
with certain choices of the two paths in the integrals. 
Then one can show that (\cite{\bf 2}, p.377)
  $$
  \wp_{12}(u)=-x_1x_2, \ \ \wp_{22}(u)=x_1+x_2. 
  \tag 1.3
  $$
If  $u$  is a point on  $\kappa^{-1}\iota(C)$,  
the  $x$-  and  $y$-coordinates of  $\iota^{-1}\kappa(u)$  
will be denoted by  $x(u)$  and  $y(u)$, respectively.  

\proclaim{\indent\sc Lemma 1.6}
If  $u$  is a point on  $\kappa^{-1}\iota(C)$,  then  
  $$
%   \frac{\wp_{11}}{\wp_{12}}(u)=
  \frac{\wp_{12}}{\wp_{22}}(u)
  =\frac{\sigma_1}{\sigma_2}(u)
  =-x(u). 
  $$
\endproclaim

\fp
These equalities are shown in  \cite{{\bf 12}, p.124} 
by using 1.4(1) and (1.3) 
(see also the proof of Proposition 2.1.1 in \cite{\bf 18}).  
By 1.2 and 1.4(1), we see the following.  

\proclaim{\indent\sc Lemma 1.7} 
If  $u\in\kappa^{-1}\iota(C)$, then 
  $$
  u_1=\frac13 {u_2}^3 + (d^{\circ}(u_2)\geq 4)
  $$
\endproclaim

\fp
Hence  $u_2$  is a local parameter at  $(0, 0)$. 
Lemma 1.8 below gives a local parameter at any point 
of  $\kappa^{-1}\iota(C)$  except  $\kappa^{-1}(O)$.

\proclaim{\indent\sc Lemma 1.8} 
Let  $v\lr{0}=(v\lr{0}_1, v\lr{0}_2)$  be a fixed point 
on  $\kappa^{-1}\iota(C)$  and 
assume  $v\lr{0}\not\in \kappa^{-1}(O)$.  
Let  $u=(u_1, u_2)$  be a variable point on  $\kappa^{-1}\iota(C)$.  
Then the variable  $u_1-v\lr{0}_1$  is a local parameter 
at  $v\lr{0}$  along  $\kappa^{-1}\iota(C)$.  
In other words, the function  $u\mapsto u_1-v\lr{0}_1$  
vanishes at  $u=v\lr{0}$  of order $1$.  
\endproclaim

{\it Proof.}  
Since
  $$
   \dfrac{d(u_2-v\lr{0}_2)}{d(u_1-v\lr{0}_1)}
  =\dfrac{du_2}{du_1}
  =\dfrac{du_2}{dx}\dfrac{dx}{du_1}=x(u), 
  $$
we have  
  $$
  u_2-v\lr{0}_2
  =x(v\lr{0})(u_1-v\lr{0}_1)+(d^{\circ}(u_1-v\lr{0}_1)\geq 2). 
  $$
Lemma 1.4(3) states that 
there exists a point  $v\lr{1}$  in  $\kappa^{-1}\iota(C)$  such that 
the function  $u\mapsto \sigma(u-v\lr{0}-v\lr{1})$  is not identically zero 
and vanishes at  $u=v\lr{0}$  of order $1$. 
Then we have
  $$
  \aligned
  &\sigma(u-v\lr{0}-v\lr{1}) \\
  &=\sigma_1(-v\lr{1})(u_1-v\lr{0}_1) + \sigma_2(-v\lr{1})(u_2-v\lr{0}_2)
   +(d^{\circ}(u_1-v\lr{0}_1, u_2-v\lr{0}_2)\geq 2)
  \endaligned
  $$
Hence the vanishing order of  $u\mapsto u_1-v\lr{0}_1$  at  $u=v\lr{0}$  
can not be higher than  $1$.  
So we see that  $u_1-v\lr{0}_1$  is a local parameter at  $v\lr{0}$.  
\qedright

\proclaim{\indent\sc Lemma 1.9}
{\rm (1)} Let  $u$  be an arbitrary point on  $\kappa^{-1}\iota(C)$.  
Then  $\sigma_2(u)$  is  $0$  if and
only if  $u$  is a lattice point, that is, 
the case  $\kappa(u)=\iota(\infty)$.  \fp
{\rm (2)} The Taylor expansion of the function  $\sigma_2(u)$  on  
$\kappa^{-1}\iota(C)$  at  $u=(0, 0)$  is of the form
  $$
  \sigma_2(u) = -{u_2}^2 + (d^{\circ}(u_2) \geq 3).
  $$
\fp
{\rm (3)} Let  $v\lr{0}=(v\lr{0}_1, v\lr{0}_2)$  be a fixed point 
on  $\kappa^{-1}\iota(C)$  and assume  $v_0\not\in \kappa^{-1}(O)$.  
The Taylor expansion of the function  $u\mapsto \sigma(u-v\lr{0})$  
on  $\kappa^{-1}\iota(C)$  at  $u=(0, 0)$  is of the form
  $$
  \sigma(u-v\lr{0})=\sigma_2(v\lr{0})u_2+(d^{\circ}(u_2)\geq 2). 
  $$
\endproclaim

{\it Proof.} If  $\sigma_2(u)=0$, 
then the second equality of 1.6 yields  $\sigma_1(u)=0$.   
This contradicts to 1.4(1), (2) and 1.1.  
So it must be  $\sigma_2(u)\neq 0$.  
The assertion (2) follows from 1.2 and 1.7.  
Because  $\sigma_2(u)$  is an even function, 
the assertion (3) follows from 1.4(1) and 1.7.  
\qedright

\proclaim{\indent\sc Definition-Proposition 1.10} 
Let  $n$  be a positive integer. If  $u\in\kappa^{-1}\iota(C)$, then
  $$
  \psi_n(u):=\frac{\sigma(nu)}{\sigma_2(u)^{n^2}}.  
  $$
is periodic with respect to  $\Lambda$. 
In other words it is a function on  $\iota(C)$. 
\endproclaim

\fp
For a proof of this, see \cite{{\bf 12}, p.124}  
or Proposition 3.2.2 in \cite{\bf 18}. 

\proclaim{\indent\sc Lemma 1.11}  
We have that  $\psi_2(u)=2y(u)$. 
\endproclaim

\fp
This is a result stated in \cite{{\bf 12}, p.128} 
(see also Lemma 3.2.4 in \cite{\bf 18}).  
We end this Section by stating the following easily shown relations 
(see Lemma 2.3.1 in \cite{\bf 18}).

\proclaim{\indent\sc Lemma 1.12}
If  $u\in \kappa^{-1}\iota(C)$  then 
  $$
  x(u)= \frac 1{{u_2}^2}+(d^{\circ}(u_2)\geq -1),\ \  
  y(u)=-\frac 1{{u_2}^5}+(d^{\circ}(u_2)\geq -4). 
  $$
\endproclaim

$ $

\par{\bf 2. A Generalization of the Formula of Frobenius-Stickelberger. } 
Let us start by stating the simplest case of our generalization. 

\proclaim{\indent\sc Proposition 2.1}
Assume  $u$  and  $v$  belong to  $\kappa^{-1}\iota(C)$.  Then 
  $$
  -\frac{\sigma(u+v)\sigma(u-v)}{\sigma_2(u)^2\sigma_2(v)^2}
  =-x(u)+x(v)
  \left(=\left|\matrix 1 &  x(u) \\
                       1 &  x(v) 
         \endmatrix\right|\right). 
  $$
\endproclaim

{\it Proof.} We give two proofs. First we use the formula of 1.5.   
For  $u\notin \kappa^{-1}\iota(C)$  and  $v\notin \kappa^{-1}\iota(C)$,  
after dividing the formula of 1.5 by  
  $$
  \wp_{22}(u)\wp_{22}(v)
=\frac{(\sigma_2(u)^2-\sigma_{22}(u)\sigma(u))}{\sigma(u)^2}\cdot
 \frac{(\sigma_2(v)^2-\sigma_{22}(v)\sigma(v))}{\sigma(v)^2}, 
  $$  
and bringing  $u$  and  $v$  close to 
any points on   $\kappa^{-1}\iota(C)$, we have 
the desired formula because of 1.4(1) and 1.6. 
\par
Our second proof is done by comparing the zeroes and poles of each side. 
If we regard the two sides as functions of  $u$, 
they are functions on  $\iota(C)$.  
We may assume  $v\not\in\kappa^{-1}(O)$.  
The zeroes of the two sides coincide and they are at  $u=v$  and  at  $u=-v$, 
and both sides have their only pole at  $u=(0, 0)$  of 
order 2 by 1.9(2) and 1.12. 
By 1.9(3), 
the coefficient of  $1/{u_2}^2$  of the left hand side is  $-1$.  
Thus we have proved the formula once again. 
\qedright

\remark{\indent\sc Remark 2.2} 
Although the formula of 1.5 is a natural generalization of 
the corresponding formula for Weierstrass' functions  $\sigma(u)$  
and  $\wp(u)$ (Example 1 in \cite{{\bf 21}, p.451}),  
the formula of 2.1 bears a striking likeness to 
the formula for elliptic functions. 
Finding this formula was the clue 
to the discovery of the formula in Theorem $2.3$. 
\endremark

Our generalization of (0.2) is the following formula.  

\proclaim{\indent\sc Theorem 2.3} 
Let  $n$  be a positive integer.  
Assume  $u(=u\lr{0})$, $u\lr{1}$, $\cdots$,  $u\lr{n}$  all 
belong to  $\kappa^{-1}\iota(C)$.  
Then 
  $$
  -\frac{\sigma(u\lr{0}+u\lr{1}+\cdots+u\lr{n})\prod_{i<j}\sigma(u\lr{i}-u\lr{j})}
        {\sigma_2(u\lr{0})^{n+1}\sigma_2(u\lr{1})^{n+1}\cdots\sigma_2(u\lr{n})^{n+1}}
  $$
is equal to
  $$
   \left|\matrix
     1  &  x(u\lr{0}) & x^2(u\lr{0}) & y(u\lr{0})            & x^3(u\lr{0}) &
          yx(u\lr{0}) & \cdots   & yx^{(n-4)/2}(u\lr{0}) & x^{(n+2)/2}(u\lr{0}) \\
     1  &  x(u\lr{1}) & x^2(u\lr{1}) & y(u\lr{1})            & x^3(u\lr{1}) &
          yx(u\lr{1}) & \cdots   & yx^{(n-4)/2}(u\lr{1}) & x^{(n+2)/2}(u\lr{1}) \\
 \vdots &  \vdots & \vdots   & \vdots            & \vdots   &
          \vdots  & \ddots   & \vdots            & \vdots           \\
     1  &  x(u\lr{n}) & x^2(u\lr{n}) & y(u\lr{n})            & x^3(u\lr{n}) & 
          yx(u\lr{n}) & \cdots   & yx^{(n-4)/2}(u\lr{n}) & x^{(n+2)/2}(u\lr{n}) \\
   \endmatrix\right|
  $$
or 
  $$
   \left|\matrix
     1  &  x(u\lr{0}) & x^2(u\lr{0}) & y(u\lr{0})            & x^3(u\lr{0}) &
          yx(u\lr{0}) & \cdots   & x^{(n+1)/2}(u\lr{0})  & yx^{(n-3)/2}(u\lr{0}) \\
     1  &  x(u\lr{1}) & x^2(u\lr{1}) & y(u\lr{1})            & x^3(u\lr{1}) &
          yx(u\lr{1}) & \cdots   & x^{(n+1)/2}(u\lr{1})  & yx^{(n-3)/2}(u\lr{1}) \\
 \vdots &  \vdots & \vdots   & \vdots            & \vdots   & 
          \vdots  & \ddots   & \vdots            & \vdots           \\
     1  &  x(u\lr{n}) & x^2(u\lr{n}) & y(u\lr{n})            & x^3(u\lr{n}) &
          yx(u\lr{n}) & \cdots   & x^{(n+1)/2}(u\lr{n})  & yx^{(n-3)/2}(u\lr{n}) \\
   \endmatrix\right|
  $$
according to whether  $n$  is even or odd.  
\endproclaim

{\it Proof.}  
Proposition 2.1 is the case of  $n=1$.  
We prove this formula by induction on  $n$.  
Suppose the points  $u\lr{0}$, $u\lr{1}$, $\cdots$, $u\lr{n}$  do 
not belong to  $\kappa^{-1}\iota(C)$.  
We know that each the function
  $$
  u\lr{j} \mapsto 
  -\frac{\sigma(u\lr{0}+u\lr{1}+\cdots+u\lr{n})\prod_{i<j}\sigma(u\lr{i}-u\lr{j})}
        {\sigma(u\lr{0})^{n+1}\sigma(u\lr{1})^{n+1}\cdots\sigma(u\lr{n})^{n+1}}
  $$
is periodic with respect to  $\Lambda$  
by 1.4(1) and the theorem of square (\cite{{\bf 17}}, Coroll.4 in p.59).  
After multiplying  $(n+1)$-st power of 
  $$
  \align
  &(-1)^{n+1}
  \frac{\wp_{22}(u\lr{0})\wp_{22}(u\lr{1})\cdots \wp_{22}(u\lr{n})}
       {\wp_{222}(u\lr{0})\wp_{222}(u\lr{1})\cdots \wp_{222}(u\lr{n})} \\
 &=\frac{(\sigma_2^2-\sigma_{22})\sigma}
       {\sigma_2^3-3\sigma\sigma_2\sigma_{22}+\sigma^2\sigma_{222}}(u\lr{0})
  \frac{(\sigma_2^2-\sigma_{22})\sigma}
       {\sigma_2^3-3\sigma\sigma_2\sigma_{22}+\sigma^2\sigma_{222}}(u\lr{1})
  \cdots
  \frac{(\sigma_2^2-\sigma_{22})\sigma}
       {\sigma_2^3-3\sigma\sigma_2\sigma_{22}+\sigma^2\sigma_{222}}(u\lr{n}), 
  \endalign
  $$
bringing  $u\lr{j}$'s close to any points of  $\kappa^{-1}\iota(C)$,  
we see by 1.4(1) the left hand side of the claimed formula is, 
as a function of each  $u\lr{j}$,  
a periodic function on  $\kappa^{-1}\iota(C)$.  
Now we regard the both sides of the claimed formula as functions 
on  $C$  of the variable  $u=u\lr{0}$,  
and regard the points  $u_j$ as points on  $C$.  
We denote the left hand side by  $f_1(u)$  
and the right hand side by  $f_2(u)$.  
We consider a divisor 
  $$
  D=(n+2)\infty-u\lr{1}-u\lr{2}-\cdots -u\lr{n}
  $$
of  $C$.  
Lemmas 1.4(2) and 1.9 show that the divisor  $(f_1)+D$  is effective, and  
1.12 shows that the divisor  $(f_2)+D$  is also effective.  
Since  $D$  is not a canonical divisor and  $C$  is of genus  $2$, 
the Riemann-Roch theorem shows 
that the space of the functions  $f$  such that  $(f)+D$  is 
effective is of dimension $1$. 
Therefore the two sides coincide up to 
a non-zero multiplicative constant.  
We know by 1.4(2), 1.9(2) and 1.9(3) that the coefficient of 
the Laurent expansion at  $u=(0,0)$  with respect to  $u_g$  
of the left hand side is just the left hand side 
of the hypothetical statement of our induction.  
We also know by 1.12 that the coefficient of 
such the Laurent expansion of the right hand side is 
just the right hand side of the hypothesis.  
Thus the formula holds.  
\qedright

\remark{\indent\sc Remark 2.4} The formula in 2.3 should be regarded 
as a generalization of the genus-one formula
  $$
  \aligned
  &\frac{\sigma(u\lr{0}+u\lr{1}+\cdots+u\lr{n})\prod_{i<j}\sigma(u\lr{i}-u\lr{j})}
        {\sigma(u\lr{0})^{n+1}\sigma(u\lr{1})^{n+1}\cdots\sigma(u\lr{n})^{n+1}} \\
  &=\left|\matrix
     1 &  x(u\lr{0}) & y(u\lr{0}) & x^2(u\lr{0}) & yx(u\lr{0}) & x^3(u\lr{0}) & \cdots  \\
     1 &  x(u\lr{1}) & y(u\lr{1}) & x^2(u\lr{2}) & yx(u\lr{2}) & x^3(u\lr{1}) & \cdots  \\
\vdots &  \vdots & \vdots   & \vdots & \vdots   & \vdots  & \ddots  \\
     1 &  x(u\lr{n}) & y(u\lr{n}) & x^2(u\lr{n}) & yx(u\lr{n}) & x^3(u\lr{n}) & \cdots  \\
   \endmatrix\right|
  \endaligned
  $$
rather than of (0.1).  
Here we suppose the sigma function  $\sigma(u)$  is associated with 
the elliptic curve defined by an equation of 
the form  $y^2=x^3+\lambda_2x^2+\lambda_1x+\lambda_0$.  
The size of the matrix above is also  $n+1$  by  $n+1$.    
This formula is easily obtained from (0.1) and explains the meaning 
of the multiplicative constant of the left hand side of (0.1).   
\endremark

$ $

\par{\bf 3. A Determinant Expression for Generalized Psi-Functions. }  
We give in this Section a generalization of the formula 
of Kiepert \cite{\bf 14}
(see also Exercise 33 in \cite{{\bf 21}, p.460} or \cite{{\bf 9}, p.186}). 
There is a pretty formula:

\proclaim{\indent\sc Lemma 3.1} Let  $u$  and  $v$  be 
belong to  $\kappa^{-1}\iota(C)$.  Then 
  $$ 
  \lim_{u_1\to v_1}\frac{\sigma(u-v)}{u_1-v_1}=1.
  $$
\endproclaim

{\it Proof.} Because of 2.1 we have
  $$
  \frac{x(u)-x(v)}{u_1-v_1}
  =\frac{\sigma(u+v)}{\sigma_2(u)^2\sigma_2(v)^2}
   \cdot
   \frac{\sigma(u-v)}{u_1-v_1}.
  $$
Now we bring  $u_1$  close to  $v_1$.  
Then the limit of the left hand side is  
  $$
  \lim_{u_1\to v_1}\frac{x(u)-x(v)}{u_1-v_1}
  =\frac{dx}{du_1}(v).
  $$
This is equal to  $2y$  by (1.2).  
The required formula follows from 1.11. 
\qedright

\remark{\indent\sc Remark 3.2}
The reader should take care to note that the left hand side in 3.1 is not    
  $$
  \lim_{h_1\to 0}\frac{\sigma(h)}{h_1} \ \
  (h \in \kappa^{-1}\iota(C)). 
  $$
\endremark

Our generalization of the formula  (0.2)  
for  $\psi_n(u)(=\sigma(nu)/\sigma_2(u)^{n^2})$  is the following.  

\proclaim{\indent\sc Theorem 3.3} 
Let  $n$  be an integer greater than $1$.  
Assume that  $u$  belongs to  $\kappa^{-1}\iota(C)$.  
Let  $j$  be  $1$  or  $2$.  
Then the following formula holds{\rm :}  
  $$
  \aligned
  &-(1!2!\cdots (n-1)!)\psi_n(u) 
% =\frac{\sigma(nu)}{\sigma_2(u)^{n^2}}
  = x^{(j-1)n(n-1)/2}(u)\times \\
  &\left|\matrix 
       x'(u)        & (x^2)'(u)         & y'(u)        & (x^3)'(u)  
  & (yx)'(u)        & (x^4)'(u)        & \cdots                    \\
       x''(u)       & (x^2)''(u)        & y''(u)       & (x^3)''(u)
  & (yx)''(u)       & (x^4)''(u)       & \cdots                    \\
       x'''(u)      & (x^2)'''(u)       & y'''(u)      & (x^3)'''(u)
  & (yx)'''(u)      & (x^4)'''(u)      & \cdots                    \\
       \vdots       & \vdots            & \vdots       & \vdots 
  & \vdots          & \vdots            & \ddots                    \\
    x\lr{n-1}(u)    & (x^2)\lr{n-1}(u)  & y\lr{n-1}(u) & (x^3)\lr{n-1}(u) 
  & (yx)\lr{n-1}(u) & (x^4)\lr{n-1}(u) & \cdots                    \\
  \endmatrix\right|.
  \endaligned
  $$
Here the size of the matrix is  $n-1$  by  $n-1$.  
The symbols  ${}'$, ${}''$, $\cdots$, ${}\lr{n-1}$  denote  
$\frac{d}{du_j}$, $\left(\frac{d}{du_j}\right)^2$, $\cdots$, 
$\left(\frac{d}{du_j}\right)^{n-1}$, respectively.     
\endproclaim

{\it Proof.} We start with replacing  $n$  by  $n-1$  in 2.3.  
If we use  $\frac{d}{du_1}=x(u)\frac{d}{du_2}$  
in the following proof for the case of  $j=1$, 
we can easily check that the case of  $j=2$  is also holds.  
We may assume that  $u\neq(0,0)$  modulo  $\Lambda$.  
For a fixed  $u$, 
we write  $h=(h_1, h_2)=u\lr{1}-u$  and varies  $u\lr{1}$.  
Then  $h_1$  is a local parameter at the point  $u$  by 1.8 
and the right hand side of 2.3 is equal to
  $$
  \aligned
  &
  \left|\matrix
    1    &      x(u)     &       x^2(u)     &   \cdots    \\
    0    &  x(u+h)-x(u)  & x^2(u+h)-x^2(u)  &   \cdots    \\
    1    &      x(u\lr{2})   &       x^2(u\lr{2})   &   \cdots    \\
  \vdots &    \vdots     &    \vdots        &   \ddots    
  \endmatrix\right|\\
  = &
  \left|\matrix
    1    &  x(u)   &  x^2(u)   & \cdots \\
    0    & x'(u)h_1+(d^{\circ}(h_1)\geq 2)
                   & (x^2)'(u)h_1+(d^{\circ}(h_1)\geq 2) 
                               & \cdots \\
    1    &  x(u\lr{2}) &  x^2(u\lr{2}) & \cdots \\
  \vdots & \vdots  &  \vdots   & \ddots
  \endmatrix\right|
  \endaligned
  $$
by Talyor's theorem in one variable. 
% We should take care to note that  $u\lr{1}$  and   $u$  are variables 
% on  $\kappa^{-1}\iota(C)$.  
After dividing both sides by  $h_1$  and taking the limit 
when  $h_1\rightarrow 0$  while keeping  $u+h\in \kappa^{-1}\iota(C)$, 
by using 3.1,  we arrive at the formula
  $$
  \aligned &
  -\frac{\sigma(2u+u\lr{2}+\cdots+u\lr{n-1})
         \sigma(u-u\lr{2})^2 \cdots \sigma(u-u\lr{n-1})^2
        \prod_{2\leq i<j}\sigma(u\lr{i}-u\lr{j})}
        {\sigma_2(u)^{2n}\sigma_2(u\lr{2})^n\cdots
         \sigma_2(u\lr{n-1})^n} \\
  &=\left|\matrix
     1     &  x(u)   &  x^2(u)   & \cdots \\
     0     &  x'(u)  & (x^2)'(u) & \cdots \\
     1     &  x(u\lr{2}) &  x^2(u\lr{2}) & \cdots \\
    \vdots &  \vdots & \vdots    & \ddots \\
   \endmatrix\right|.
  \endaligned
  $$
Applying Talyor's theorem in one variable again we see that 
the right hand side with  $h=u\lr{2}-u$  is given by 
  $$
  \aligned
  &
  \left|\matrix
     1   &  x(u)       &  x^2(u)        & \cdots \\
     0   &  x'(u)      & (x^2)'(u)      & \cdots \\
     0   & x(u+h)-\left(x(u)+x'(u)h_1\right)
                       & x^2(u+h)-\left(x^2(u)+(x^2)'(u)h_1\right) 
                                        & \cdots \\
     1   &  x(u\lr{3})     &  x^2(u\lr{3})      & \cdots \\
  \vdots &  \vdots     &  \vdots        & \ddots \\
  \endmatrix\right| \\
  =&\left|\matrix
     1   & x(u)   &  x^2(u)     & \cdots \\
     0   & x'(u)  & (x^2)'(u)   & \cdots \\
     0   & \frac 1{2!}x''(u)(h_1)^2+(d^{\circ}(h_1)\geq 3) 
                  & \frac 1{2!}(x^2)''(u)(h_1)^2+(d^{\circ}(h_1)\geq 3)
                                & \cdots \\
     1   & x(u\lr{3}) &  x^2(u\lr{3})   & \cdots \\
  \vdots & \vdots &  \vdots     & \ddots \\
   \endmatrix\right|.
  \endaligned
  $$ 
After dividing both sides by  $(h_1)^2$  and taking the limit 
when  $h\lr{1}\rightarrow 0$,  
by applying 3.1 we have 
  $$
  \aligned &
  -(1!2!)\frac{\sigma(3u+u\lr{3}+\cdots+u\lr{n-1})
         \sigma(u-u\lr{3})^3 \cdots \sigma(u-u\lr{n-1})^3
         \prod_{3\leq i<j}\sigma(u\lr{i}-u\lr{j})}
        {\sigma_2(u)^{3n}\sigma_2(u\lr{3})^n\cdots
         \sigma_2(u\lr{n-1})^n} \\
  &=\left|\matrix
     1     &  x(u)     &  x^2(u)    & \cdots \\
     0     &  x'(u)    & (x^2)'(u)  & \cdots \\
     0     &  x''(u)   & (x^2)''(u) & \cdots \\
     1     &  x(u\lr{3})   &  x^2(u\lr{3})  & \cdots \\
    \vdots &  \vdots   & \vdots     & \ddots \\
   \endmatrix\right|.
  \endaligned
  $$
Using similar operations repeatedly, we have the formula for  $\psi_n(u)$. 
\qedright

\remark{\indent\sc Example} 
We have
  $$
  \align
  \allowmathbreak
  x'(u)&=\frac{d}{\ du_1}x(u)
        =1\Big/\frac{\ du_1}{dx}(u)=\frac{2y}{x}(u),\ \ 
  y'(u) =\frac {Df(x)}x, \\
  x''(u)&=\frac d{\ du_1}\left(\frac{2y}{x}\right)
         =2\frac{x\cdot Df(x)-2y^2}{x^3}, \  \
 (x^2)'(u)=4y(u), \  \  
 (x^2)''(u)=\frac{4}xDf(x), \ \  
  \endalign
  $$
where  $f(x)$  is as in the beginning of Section 1 and 
$Df(x)=\frac{d}{dx}f(x)$.  
Then  $-\psi_3(u)=-\sigma(3u)/{\sigma_2(u)^9}=8y^3(u)$  by 3.3. 
This example is mentioned in the proof of Lemma 2 (e) of \cite{\bf 12}, too, 
and the proof here is different from the one there 
(see also Lemma 3.2(2) in \cite{\bf 18}).  

\endremark

$ $

$ $

%  \centerline{---------------------------}

\enddocument
\bye